\documentclass{article}

 \usepackage{graphicx}

\usepackage{multirow} 
\usepackage{amsmath, amssymb, amsthm, enumerate}
\usepackage{amsmath,amssymb,amsthm,enumerate,mathrsfs}
\newtheorem{theorem}{Theorem}[section]

\newtheorem{proposition}[theorem]{Proposition}
\theoremstyle{definition}

\usepackage{enumitem}
\usepackage{threeparttable}
\usepackage{amsthm}
\usepackage{bm} 

\usepackage{algorithm}%
\usepackage{algorithmicx}
\usepackage{algpseudocode}%
\usepackage{listings}%
\usepackage{xcolor}%
\usepackage{textcomp}%
\usepackage{manyfoot}%
\usepackage{booktabs}%

\numberwithin{equation}{section}
\usepackage{caption}
\begin{document}
\makeatletter

\begin{center}
\large{\bf An Efficient Solution Method for Solving Convex Separable Quadratic Optimization Problems
}\\[5mm]
\large{Dedicated to Professor Terry Rockafellar on the occasion of his 90th birthday}
\end{center}\vspace{5mm}
\begin{center}
\textsc{Shaoze Li, Junhao Wu, Cheng Lu, Zhibin Deng, Shu-Cherng Fang}\end{center}

\vspace{2mm}

{\footnotesize{
\noindent\begin{minipage}{14cm}
{\bf Abstract:}
Convex separable quadratic optimization problems occur in many practical applications. In this paper, based on an iterative resolution scheme of the KKT system, we develop an efficient method for solving a quadratic programming problem with a convex separable objective function subject to multiple convex separable constraints. We show that the proposed approach leads to a dual coordinate ascent algorithm and provide a convergence proof. Numerical experiments support the superior performance of the proposed method to that of the Gurobi solver, especially for solving large-size convex separate quadratic programming problems.
\end{minipage}
 \\[5mm]


\hbox to14cm{\hrulefill}\par}}


\section{Introduction}\label{sec1}
In this paper, we consider the following convex separable quadratic  program with multiple separable quadratic constraints (SQPQC):
\begin{align*}
\min _{\bm{y} \in \mathbb{R}^n} f(\bm{y}) & =\bm{y}^T \Delta \bm{y}+\bm{\alpha}^T \bm{y} \\
\text { s.t. } g_{i}(\bm{y})& = \bm{y}^T \Theta_i \bm{y} +\bm{\beta}_i ^T \bm{y} +\sigma_i \leq 0,~ i=1,\dots,m,   \label{SQPQC}\tag{SQPQC}\\ 
\bm{y} & \in[\bm{l}, \bm{u}],
\end{align*}
where $\Delta$ and $\Theta_i$ are $n \times n$ diagonal matrices, i.e., $\Delta=\operatorname{diag}\left(\delta_1, \ldots, \delta_n\right)$, $\Theta_i=\operatorname{diag}\left(\theta_i^1, \ldots, \theta_i^n\right)$, $\bm{\alpha}=\left(\alpha_1, \ldots, \alpha_n\right)^{T} \in \mathbb{R}^n$, $\bm{\beta}_i=\left(\beta_i^1, \ldots, \beta_i^n\right)^{T} \in \mathbb{R}^n$, $\sigma_i \in \mathbb{R}$, $\bm{l}=\left(l_1, \ldots, l_n\right)^{T} \in \mathbb{R}^n$ and $\bm{u}=\left(u_1, \ldots, u_n\right)^{T} \in \mathbb{R}^n$. Throughout the paper, we assume that the following two conditions  hold:
\begin{itemize}
\setlength{\itemindent}{-5mm}
\item[$\bullet$] Assumption 1: $f(\bm{y})$ is strongly convex and $g_{i}(\bm{y})$ is  convex, i.e., $\delta_{j}>0$ and $\theta^j_{i}\geq 0$ for $i=1,\dots,m$ and $ j=1,\dots,n$. \label{assump1}

\item[$\bullet$] Assumption 2: There exists  $\hat{\bm y} \in (\bm{l},\bm{u})$ such that $g_{i}(\hat{\bm{y}})<0~\text{for} ~ i=1,\dots,m.$\label{assump2}
\end{itemize}
Assumption 1 assures that the objective value is bounded from below over the feasible region.  Assumption 2 says that Slater's condition holds for problem \eqref{SQPQC}. 


 Problem \eqref{SQPQC} has various practical applications, such as the resource allocation problems \cite{bitran1981disaggregation,hochbaum1995strongly} and multi-commodity network flows \cite{shetty1990parallel}. Some of these applications contain numerous decision variables.  It is known that problem \eqref{SQPQC} can be cast into a second-order cone programming problem \cite{fukushima2003sequential, lobo1998applications} and then solved using interior point methods \cite{fukushima2003sequential, lobo1998applications, mehrotra1991method, monteiro2000polynomial,  tsuchiya1999convergence}. It can also be solved by the sequential quadratic programming algorithm \cite{pjo1,pjo2}. However, when solving large-size problems, the computational burden of these methods remains a critical issue. This motivates us to develop an efficient algorithm tailored for solving  large-size problem \eqref{SQPQC}.

There are many works related to problem \eqref{SQPQC} due to its broad applicability.
When problem \eqref{SQPQC} degrades to the problem with only linear constraints and boxed variables, that is, $\Theta_{i}=0$, $i=1,...,m$, it becomes a linearly constrained quadratic separable programming problem. Megiddo and Tamir \cite{megiddo1993linear} showed that the Lagrangian dual method derived from the multidimensional search procedure proposed by Megiddo \cite{megiddo1984linear} yields a linear-time algorithm for this subclass of problem \eqref{SQPQC}. Nevertheless, this search method, which relies on the linear property of the constraints, fails for problems with quadratic constraints. Tseng \cite{tseng1990dual} presented a dual ascent method for this special problem and provided the convergence property of the proposed method under several assumptions. However, the proof for problems with linear constraints cannot be directly extended to problems with quadratic ones.  When problem \eqref{SQPQC} is subjected to only one linear constraint and boxed variables,  Dai and Fletcher \cite{dai2006new} proposed a tailored secant method with good practical performance to solve  this problem. Then Cominetti et al. \cite{cominetti2014newton} further developed a semismooth Newton method to solve this special problem with better performance than that in \cite{dai2006new}.

 Additional research works have extended the study of problem \eqref{SQPQC}  to the nonconvex case.  When the problem with a separable nonconvex quadratic objective function has only one linear constraint and boxed variables, Edirisinghe and Jeong \cite{edirisinghe2016efficient,edirisinghe2017tight} presented fast linear-time complexity procedures to compute tight lower and upper bounds for  large-scale instances.  Li et al. \cite{li2023efficient} further proposed an efficient global algorithm based on the KKT branching to solve large-scale nonconvex quadratic knapsack problems.   Furthermore,  when the problem with  nonconvex separable quadratic objective function has only one nonconvex quadratic constraint with boxed variables, Chen et al. \cite{Chen2014} introduced a lower-bound algorithm using Lagrangian relaxation, while Luo et al. \cite{Luo2023} further developed an efficient global algorithm in a branch-and-bound framework.

In this paper,  we start from the dual solution method for problem \eqref{SQPQC} with one single quadratic constraint to exploit some distinct properties of this subclass,  such as the differentiability of the dual function of problem \eqref{SQPQC} and the uniqueness of the optimal solution of the dual problem. Leveraging on the properties, we then develop a highly efficient algorithm for solving  problem \eqref{SQPQC} with multiple quadratic constraints and prove the convergence of the proposed algorithm.

The paper is organized as follows. In Sect. \ref{sec2}, we introduce an efficient algorithm for solving problem \eqref{SQPQC} with a single separable convex quadratic constraint. Subsequently, in Sect. \ref{sec3}, we develop an iterative method to address  problem \eqref{SQPQC} with multiple constraints. The convergence of the proposed algorithm is  established in Sect. \ref{sec3}. Computational experiments  in Sect. \ref{sec4}  demonstrate the promising performance of the proposed algorithm.

\section { Convex quadratic optimization with a single separable convex quadratic constraint}\label{sec2}

In this section, we  consider a separable convex quadratic program with only  one quadratic constraint and  boxed variables in  the following form:
\begin{align*}
\min _{\bm{y} \in \mathbb{R}^n} f(\bm{y}) & =\bm{y}^T \Delta \bm{y}+\bm{\alpha}^T \bm{y}  \\ \tag{SQPQC1}\label{SQPQC1} 
\text { s.t. } g(\bm{y})& = \bm{y}^T \Theta \bm{y} +\bm{\beta} ^T \bm{y} +\sigma \leq 0,  \\
\bm{y} & \in[\bm{l}, \bm{u}],
\end{align*}
where $\Theta$ is a diagonal positive semidefinite matrix.
We first introduce the existing results for problem \eqref{SQPQC1} in \cite{Chen2014, Luo2023} in Subsection \ref{2.1}. Then, based on the existing results, we derived some new properties and the solution algorithm for problem \eqref{SQPQC1} in Subsection \ref{2.2}.

\subsection{  The existing results for problem \eqref{SQPQC1}.}\label{2.1}
In this subsection, we briefly review an important result of the dual problem of \eqref{SQPQC1} in \cite{Chen2014, Luo2023}, which will be used to develop the solution algorithm for problem \eqref{SQPQC1} in Subsection \ref{2.2}. 

Define the dual function of \eqref{SQPQC1}  as follows:
\begin{align*}
D(\lambda)=&\min _{\bm{y} \in[\bm{l}, \bm{u}]} \{f(\bm{y})+\lambda g(\bm{y})\}  \\
=&\min _{\bm{y} \in[\bm{l}, \bm{u}]}\{ \bm{y}^T \Delta \bm{y}+\bm{\alpha}^T \bm{y}+\lambda (\bm{y}^T \Theta \bm{y} +\bm{\beta}^T \bm{y} +\sigma)\}\\
=&\min _{\bm{y} \in[\bm{l}, \bm{u}]} \{\bm{y}^T(\Delta+\lambda\Theta)\bm{y}+(\bm{\alpha}+\lambda\bm{\beta})^T \bm{y}+\lambda\sigma\},
\end{align*}
where $\lambda$ is  the dual multiplier of the constraint $g(\bm{y})\leq 0$. Then, the dual problem of \eqref{SQPQC1} is expressed as
\begin{equation*}\label{DSQPQC1}\tag{DSQPQC1}
    \max_{\lambda\geq 0} ~D(\lambda). 
\end{equation*}
For any given $\lambda$,  define $\bm{y}^*(\lambda)$ to be  the minimum solution to the problem of $\text{min}_{\bm{y} \in[\bm{l}, \bm{u}]} \left\{f(\bm{y})+\lambda g(\bm{y})\right\}$. In fact, the closed-form expression for $\bm{y}^*(\lambda)$ is available in the following proposition. 
\begin{proposition}(Proposition 2 in \cite{Luo2023})\label{lemma1} The closed-form expression for $\bm{y}^*(\lambda)=\{y_{1}^*(\lambda),\dots,y_{n}^*(\lambda)\} $ is given by $${y}_{j}^*(\lambda)=\max\left\{l_{j},~\min\left\{u_{j}, -\frac{\alpha_{j}+\lambda\beta_{j}}{2(\delta_{j}+\lambda\theta_{j})} \right\} \right\}\text{ for } j=1,\dots,n.$$
\end{proposition}


\subsection{The solution algorithm for problem \eqref{SQPQC1}.}\label{2.2}
In this subsection, we first investigate the properties of the dual function $D(\lambda)$ under Assumptions 1 and 2. These properties are essential for the proof of the convergence of the proposed algorithm in Sect. \ref{sec3}. Then we develop a solution algorithm for problem \eqref{SQPQC1} following the methodology proposed in \cite{Luo2023}.

\begin{proposition}\label{prop2.3}
Under Assumption 1, $D(\lambda)$ is continuously differentiable and $D'(\lambda)=g(\bm{y}^{*}(\lambda))$.
\end{proposition} 
\begin{proof}
  Based on Assumption 1 and Proposition \ref{lemma1}, for any given $\lambda\ge 0$, the optimal solution $\bm{y}^{*}(\lambda)$ is unique. Since $f(\bm{y})$ and $g(\bm{y})$ are continuous convex functions, it follows from  Theorem 35.8 in \cite{rockafellar1997convex} that $D(\lambda)$ is continuously differentiable and $D'(\lambda)=g(\bm{y}^{*}(\lambda))$.
\end{proof}


{\begin{proposition}\label{Prop2.4}
 Under Assumptions 1 and 2, if the  condition  
\begin{equation*}\tag{C1}\label{C1}
  \left\{\bm{{y}} \mid g(\bm{{y}})=0\right\} \bigcap \left\{\bm{{y}} \mid {y}_{j}\in \left\{l_{j},u_{j},-\frac{\beta_{j}}{2\theta_{j}}\right\},~ j=1,\dots, n\right\}=\emptyset  
\end{equation*}
holds for problem \eqref{SQPQC1}, then $g(\bm{y}^{*}(\lambda))$ is a nonincreasing function with respect to $\lambda$ and the dual problem \eqref{DSQPQC1} has a unique optimal solution $\lambda^*$. 
Furthermore, if $g(\bm{y}^*(0))< 0$, the optimal solution $\lambda^*$ is equal to $0$. If $g(\bm{y}^*(0))\ge 0$, then $\lambda^*$  is the solution of the equation $g(\bm{y}^*(\lambda))=0$ over $\lambda\ge 0$. 
\end{proposition} 

\begin{proof}
To prove the monotonicity of $g(y^*(\lambda))$, consider two different parameters $\lambda'$ and $\lambda''$. By the optimality of $\bm{y}^{*}(\lambda')$  and $\bm{y}^{*}(\lambda'')$, we have
\begin{align*}
  f(\bm{y}^{*}(\lambda'))+\lambda' g(\bm{y}^{*}(\lambda'))\leq f(\bm{y}^{*}(\lambda''))+\lambda' g(\bm{y}^{*}(\lambda'')),\\
f(\bm{y}^{*}(\lambda''))+\lambda'' g(\bm{y}^{*}(\lambda''))\leq f(\bm{y}^{*}(\lambda'))+\lambda'' g(\bm{y}^{*}(\lambda')).  
\end{align*}
Adding the above two inequalities, we obtain
\begin{equation*}
    (g(\bm{y}^{*}(\lambda'))-g(\bm{y}^{*}(\lambda'')))(\lambda'-\lambda'')\leq 0.
\end{equation*}
This proves that $g(\bm{y}^{*}(\lambda))$ is a nonincreasing function in $\lambda$.

Now we prove the uniqueness of the dual optimal solution. Since problem \eqref{SQPQC1}  is a convex problem satisfying Slater's condition (Assumption 2) and has a bounded function value (Assumption 1), the strong duality between \eqref{SQPQC1} and \eqref{DSQPQC1} holds, and the optimal solution of \eqref{DSQPQC1} exists. 

In the case of $g(\bm{y}^{*}(0))<0$, by the monotonicity of function $g(\bm{y}^*(\lambda))$, we have $g(\bm{y}^{*}(\lambda))<0$ for all $ \lambda\geq 0 $. Since $D'(\lambda)=g(\bm{y}^{*}(\lambda))$ due to Prop. \ref{prop2.3}, we deduce that the dual function is strictly decreasing in this case. Hence, the unique optimal solution to problem \eqref{DSQPQC1} is $\lambda^*=0$.  

If $g(\bm{y}^{*}(0))\ge 0$, we  show the uniqueness of the solution to $D'(\lambda) = g(\bm{y}^{*}(\lambda))=0$ over $\lambda\geq 0$, which in turn implies that problem \eqref{DSQPQC1} possesses a unique optimal solution.  Denote the optimal solution of problem  \eqref{DSQPQC1} as $\lambda^*$.  By the optimality condition,  $g(\bm{y}^{*}(\lambda))=0$ holds at $\lambda^*$. Condition \eqref{C1} assures that there exists at least one index $k \in \{1,\dots, n\}$ such that $y^*_{k}(\lambda^*)\notin \{l_{k},u_{k},-\beta_{k}/(2\theta_{k})\}$.  Recalling that ${y}_{k}^*(\lambda)=\max\{l_{k},~\min\{u_{k}, -\frac{\alpha_{k}+\lambda\beta_{k}}{2(\delta_{k}+\lambda\theta_{k})}\}\}$ in Prop. \ref{lemma1}, we can deduce that $y^*_{k}(\lambda^*)= -\frac{\alpha_{k}+\lambda^*\beta_{k}}{2(\delta_{k}+\lambda^*\theta_{k})}$.  If ${\alpha_{k}}/{\beta_{k}}= {\delta_{k}}/{\theta_{k}}$, then we have $y^*_{k}(\lambda^*)= -\frac{\alpha_{k}+\lambda^*\beta_{k}}{2(\delta_{k}+\lambda^*\theta_{k})}=-{\beta_{k}}/2{\theta_{k}}=-{\alpha_{k}}/2{\delta_{k}}$, which contradicts the condition \eqref{C1}. Hence, we conclude that $y^*_{k}(\lambda^*)= -\frac{\alpha_{k}+\lambda^*\beta_{k}}{2(\delta_{k}+\lambda^*\theta_{k})}$ and ${\alpha_{k}}/{\beta_{k}}\neq {\delta_{k}}/{\theta_{k}}$. Note that $-\frac{\alpha_{k}+\lambda\beta_{k}}{2(\delta_{k}+\lambda \theta_{k})}$ is a decreasing linear function for $\theta_{k}=0$, and $ -\frac{\alpha_{k}+\lambda\beta_{k}}{2(\delta_{k}+\lambda \theta_{k})}= -\frac{\beta_{k}}{\theta_{k}}+\frac{\beta_{k}/{\theta_{k}}-{\alpha_{k}}/{\delta_{k}} }{2(\delta_{k}+\lambda \theta_{k})\delta_{k}}$ is a strictly monotone function with $\lambda$ for $\theta_{k}\neq 0$, implying that $y^*_{k}(\lambda^*)\neq y^*_{k}(\lambda)$ and $\bm{y}^*(\lambda^*)\neq \bm{y}^*(\lambda)$ hold for all $\lambda \neq \lambda^*$. 
  Furthermore, since $f(\cdot)+\lambda g(\cdot) $ is a strictly convex function, we have the following strict inequalities for  $\lambda^*$  and $\lambda, (\lambda\neq \lambda^*$):
\begin{align*}
  f(\bm{y}^{*}(\lambda^*))+\lambda^* g(\bm{y}^{*}(\lambda^*))< f(\bm{y}^{*}(\lambda))+\lambda^* g(\bm{y}^{*}(\lambda)),\\
f(\bm{y}^{*}(\lambda))+\lambda g(\bm{y}^{*}(\lambda))<f(\bm{y}^{*}(\lambda^*))+\lambda g(\bm{y}^{*}(\lambda^*)).  
\end{align*}
 Adding the above two inequalities, we have 
\begin{align*}
    &~(g(\bm{y}^{*}(\lambda^*))-g(\bm{y}^{*}(\lambda)))(\lambda^*-\lambda)< 0\\
  \Leftrightarrow ~&~ -g(\bm{y}^{*}(\lambda))(\lambda^*-\lambda)< 0.
\end{align*}
Then we obtain that $g(\bm{y}^{*}(\lambda))$ strictly decreases in $\lambda^*$ and $g(\bm{y}^{*}(\lambda))\neq 0$ for $ \lambda\neq \lambda^*$. Therefore, the dual problem \eqref{DSQPQC1} has a unique optimal solution, and $\lambda^*$  is the solution of system $g(\bm{y}^*(\lambda))=0,~\lambda\ge 0$.   This completes the proof. 
\end{proof} 

{ \color{red} In Proposition 2.3, we introduce the additional condition \eqref{C1} for the uniqueness of the optimal solution to problem \eqref{DSQPQC1}. Actually, this condition is required to guarantee the convergence of our proposed algorithm for solving \eqref{SQPQC} in the next section. Intuitively,  Condition \eqref{C1} plays a similar role as constraint qualifications. Here, we use an example to demonstrate. Consider the following set of constraints:  
\begin{align*}\tag{SC} \label{SC}
& g(\bm{y})=(y_1-1)^2+y_2^2-1\le0,\\
&\frac{1}{2} \le y_1 \le \frac{3}{2},-\frac{\sqrt{3}}{2} \le y_2 \le 1. 
\end{align*}
 \begin{figure}[h!]
    \centering
    \includegraphics[width=0.4\textwidth]{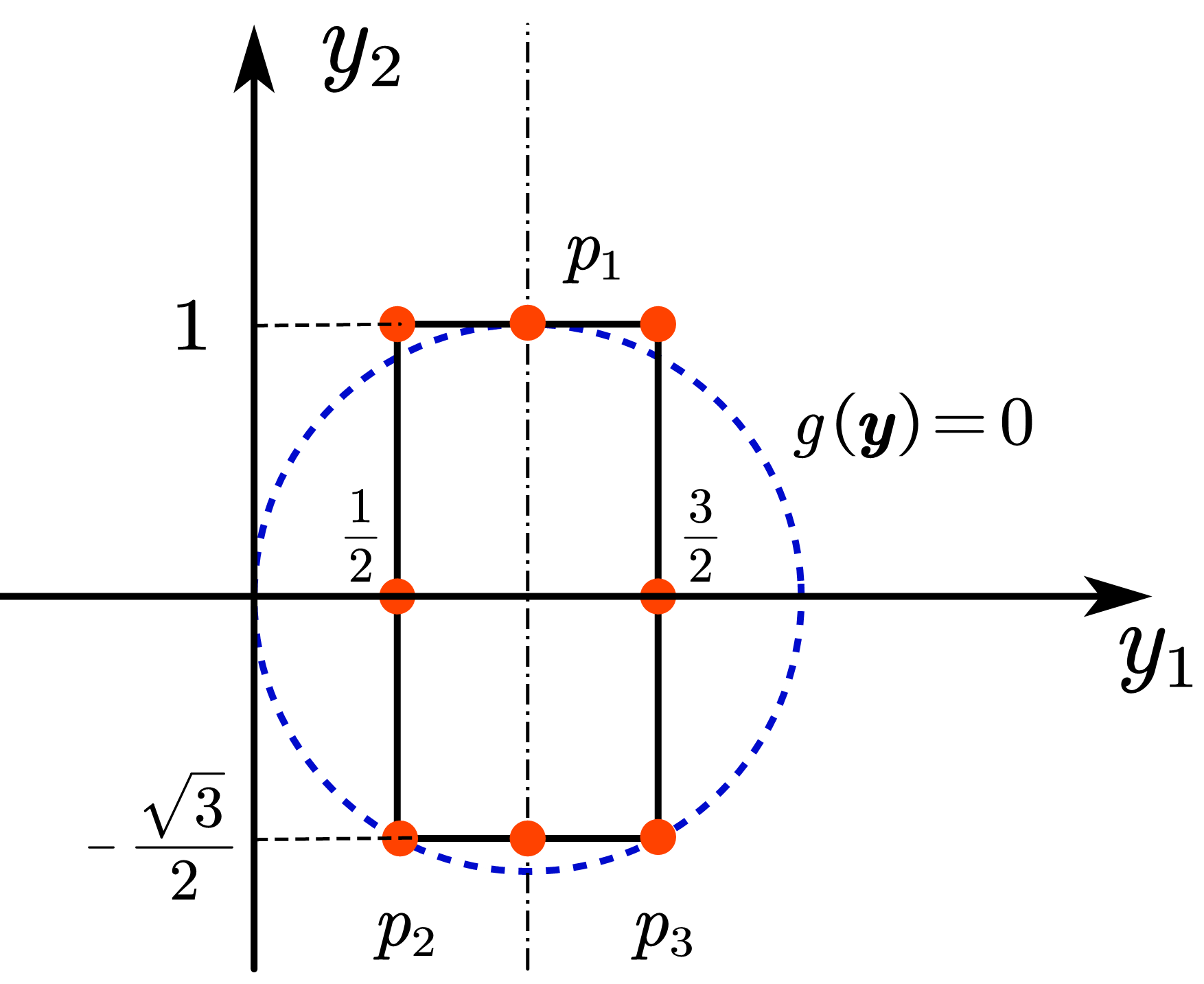}
    \caption{The boundaries of constraints in \eqref{SC}.}
    \label{fig:1}
\end{figure}

The boundaries of $g(\bm{y})\le 0$ and boxed constraints are shown in Figure \ref{fig:1}.  We indicate  $g(\bm{y})= 0$ with dashed lines, and the boxed constraints with solid lines in Figure \ref{fig:1}.   
 The dots are the elements in the set $\left\{\bm{{y}} \mid {y}_{j}\in \{l_{j},u_{j},-\frac{\beta_{j}}{2\theta_{j}}\}\right.$, $\left. j=1,\dots, n\right\}$.
In this case, we have $$\left\{\bm{{y}} \mid g(\bm{{y}})=0\right\} \bigcap \left\{\bm{{y}} \mid {y}_{j}\in \left\{l_{j},u_{j},-\frac{\beta_{j}}{2\theta_{j}}\right\},~ j=1,\dots, n\right\}=\{p_1,p_2,p_3\}.$$
Geometrically, at point $p_1$,  the gradients of active constraints $g(\bm{{y}}) \le 0$ and $y_{2} \le 1$ are collinear, implying that the linear independence constraint qualification is not satisfied. Similarly, at point $p_2$ (or $p_3$), there are three active constraints, and the linear independence constraint qualification does not hold, either. In practice, this condition holds with a high probability since the measure of the set in condition (C1) is zero.  }

Now we have established the strong duality between \eqref{SQPQC1} and \eqref{DSQPQC1}. Furthermore, leveraging Prop. \ref{Prop2.4}, we can design a solution algorithm for problem \eqref{SQPQC1}.  The intuition of our algorithm for \eqref{SQPQC1} is similar to Algorithm ILBSSA in \cite{Luo2023}, which is designed for computing a lower bound of non-convex quadratic optimization problems with a single quadratic constraint and boxed variables. We first solve the dual problem \eqref{DSQPQC1} according to the first-order optimal condition. Specifically, we  check the value of  $g(\bm{y}^*(0))$. If $g(\bm{y}^*(0))< 0$, then we set $\lambda^*=0$.  Otherwise,  we  search the solution $\lambda^*$ of the nonlinear equation $g(\bm{y}^*(\lambda))= 0$ over $\lambda\ge 0$ using the bisection method. After obtaining $\lambda^*$, we then recover the optimal solution $y^*$ of primal problem \eqref{SQPQC1} from $\lambda^*$ based on  Prop. \ref{lemma1}. The algorithm scheme for solving problem \eqref{SQPQC1} is presented in Algorithm \ref{alg1}.}

\begin{algorithm}[htbp!]
{\caption{for  \eqref{SQPQC1} }\label{alg1}
\begin{algorithmic}[1] 
\State \textbf{Input}: Problem parameters: $\Delta,~\bm{\alpha},~\Theta,~\bm{\beta},~\sigma,~\bm{l},~\bm{u},~n,~\epsilon$ and $k=0$.
\State Compute ${y}_{j}^*(0)=\max\left\{l_{j},~\min\left\{u_{j}, -\frac{\alpha_{j}}{2\delta_{j}} \right\} \right\}\text{ for } j=1,\dots,n.$ and  $g(\bm{y}^*(0))$.
\If{$g(\bm{y}^*(0))<0$}
\State $\lambda^*=0$.
\Else 
\State Set $\overline{\lambda}$=0 and find a large $\underline{\lambda}$ such that $g(\bm{y}^*(\underline{\lambda}))<0$.
\State Set ${\lambda}_m=\frac{\underline{\lambda}+\overline{\lambda}}{2}$.

\While{$|g(\bm{y}^*({\lambda}_m))|>\epsilon$}

\If{$g(\bm{y}^*(\lambda_m))>0$}
\State{$\overline{\lambda}=\lambda_{m}$},
\ElsIf{$g(\bm{y}^*(\lambda_m))<0$}
\State{$\underline{\lambda}=\lambda_{m}$},
\EndIf
\State ${\lambda}_m=\frac{\underline{\lambda}+\overline{\lambda}}{2}$,
\EndWhile

\EndIf

\While{$k \leq n$}
    \If{$-\frac{\bm{\alpha}_{k}+\lambda^{*}\bm{\beta}_{k}}{2(\delta_{k}+\lambda^{*}\theta_{k})}\geq u_{k}$}
    \State{$y_k^*=u_{k}$},
    \ElsIf{$-\frac{\bm{\alpha}_{k}+\lambda^{*}\bm{\beta}_{k}}{2(\delta_{k}+\lambda^{*}\theta_{k})}\leq l_{k}$}
    \State{$y_k^*=l_{k}$},
    \ElsIf{$l_{k}\leq -\frac{\bm{\alpha}_{k}+\lambda^{*}\bm{\beta}_{k}}{2(\delta_{k}+\lambda^{*}\theta_{k})}\leq u_{k}$}
            \State{$y_k^*=-\frac{\bm{\alpha}_{k}+\lambda^{*}\bm{\beta}_{k}}{2(\delta_{k}+\lambda^{*}\theta_{k})}$},
    \EndIf

\EndWhile
\State \textbf{Output}:  The optimal dual multiplier $\lambda^{*}$ and the optimal solution $\bm{y}^*(\lambda^{*})$. 
\end{algorithmic}}
\end{algorithm}


\section{Convex quadratic optimization with separable convex quadratic constraints} \label{sec3}
In this section, we extend  Algorithm \ref{alg1} to design an efficient algorithm for solving convex quadratic programs with multiple quadratic constraints. Notice that the algorithms for solving single constrained problems in \cite{Chen2014} and \cite{Luo2023} cannot be extended directly to our problem. The limitation comes from the fact that they solve the problem from the dual side and search for the optimal dual multiplier in the one dimensional dual space.  This approach becomes inefficient as more constraints involved inducing a higher dimentional dual space search. In this section, we propose a dual coordinate ascent algorithm based on the KKT conditions to solve the problem \eqref{SQPQC} with multiple constraints.

Consider the following KKT conditions of problem \eqref{SQPQC}:
\begin{align*}\label{KKT}
    &2(\Delta+\sum_{i=1}^{m}\lambda_{i}\Theta_{i})\bm{y}+\bm{\alpha}+\sum_{i=1}^{m}\lambda_{i}\bm{\beta_{i}}-\bm{\underline{\eta}}+\bm{\overline{\eta}} =0,\\
    &\lambda_{i}g_{i}(\bm{y})=0~,~~~~ i=1,\dots,m ,\\
    &\underline{\eta}_{j} (y_{j}-l_{j})=0,~\underline{\eta}_{j} \geq 0, ~~~~~ j=1,\dots,n,\\\tag{KKT} 
    & \overline{\eta}_{j} (y_{j}-u_{j})=0,~\overline{\eta}_{j} \geq 0,~~~~ j=1,\dots,n,\\
    &\lambda_{i}\geq 0~,~~~~ i=1,\dots,m,\\
    &g_{i}(\bm{y}) = \bm{y}^T \Theta_{i} \bm{y} +\bm{\beta}_{i} ^T \bm{y} +\sigma_{i} \leq 0, ~~~~ i=1,\dots,m, \\
&\bm{y}  \in[\bm{l}, \bm{u}],  
\end{align*}
where $\lambda_{i}$ is the dual multiplier of  $g_{i}(\bm{y})$, and $\underline{\eta}_{j}$ and $\overline{\eta}_{j}$ are the dual multipliers of $y_{j}\geq l_{j}$ and $y_{j}\leq u_{j}$, respectively. According to Assumption 2,  Slater's condition holds for \eqref{SQPQC}, and the solution of the KKT system is the solution of \eqref{SQPQC}. 

To efficiently solve the KKT condition, we consider applying an iterative scheme. In each iteration,  we focus on the subsystem that contains only one dual variable $\lambda_{k}$ by fixing $\lambda_{i}, ~i \in list(k) := \{1,\dots,m\}\backslash k$ as constants.  Specifically, in each iteration, we solve the following subsystem:
\begin{align*}\label{SKKT}
    &2(\Delta+\sum_{i\in list(k)}\hat{\lambda}_{i}\Theta_{i}+{\lambda}_{k}\Theta_{k})\bm{y}+\bm{\alpha}+\sum_{i\in list(k) }\hat{\lambda}_{i}\bm{\beta}_{i}+{\lambda}_{k}\bm{\beta}_{k}-\bm{\underline{\eta}}+\bm{\overline{\eta}} =0,\\
    &\lambda_{k}g_{k}(\bm{y})=0~,\\
    &\underline{\eta}_{j} (y_{j}-l_{j})=0,~\underline{\eta}_{j} \geq 0, ~~~~~ j=1,\dots,n,\\\tag{SKKT} 
    & \overline{\eta}_{j} (y_{j}-u_{j})=0,~\overline{\eta}_{j} \geq 0,~~~~ j=1,\dots,n,\\
    &\lambda_{k}\geq 0~,\\
    &g_{k}(\bm{y}) = \bm{y}^T \Theta_{k} \bm{y} +\bm{\beta}_{k} ^T \bm{y} +\sigma_{k} \leq 0, \\
&\bm{y}  \in[\bm{l}, \bm{u}],  
\end{align*}
where $\hat{\lambda}_{i}\geq 0,~ i \in list(k)$, are fixed constants. Notice that the subsystem \eqref{SKKT} is the KKT condition of the following problem \eqref{SSQPQC}:
\begin{align*}
   \min _{\bm{y} \in \mathbb{R}^n} ~&~ f(\bm{y})+\sum_{i\in list(k)} \hat{\lambda}_{i}g_{i}(\bm{y}) \\
\text {s.t.} ~&~ g_{k}(\bm{y}) = \bm{y}^T \Theta_k \bm{y} +\bm{\beta}_k ^T \bm{y} +\sigma_k \leq 0, \label{SSQPQC}\tag{SSQPQC}\\
~&~ \bm{y}\in[\bm{l}, \bm{u}].
\end{align*}
Under Assumptions 1 and 2, problem \eqref{SSQPQC} is strongly convex and Slater's condition holds, hence its optimal solution and the corresponding dual optimal solution must exist and solve the subsystem \eqref{SKKT}. We point out that problem \eqref{SSQPQC} can be efficiently solved by Algorithm \ref{alg1} in Sec. \ref{sec2}, yielding an optimal dual multiplier $\lambda^*$ and the optimal solution $\bm{y}^*(\lambda^*)$.




Now, we propose a solution method for problem \eqref{SQPQC}. The algorithmic procedures are presented in the Main Algorithm, in which $K$ is the maximal number of iterations, $Iter$ is the index of iteration number. In the proposed algorithm, the multiplier $\bm{\lambda}$ is updated in a component-wise cyclic order. The algorithm continues until the dual multiplier $\lambda_{k}$ and the corresponding solution $\bm{y}^*(\lambda_{k})$ in current iteration satisfy system \eqref{KKT} within a given tolerance $\epsilon$. 

\begin{algorithm}[H]
\captionsetup{labelformat=empty} 
\caption{\textbf{Main Algorithm}  for  \eqref{SQPQC} }\label{alg3}
\begin{algorithmic}[1]
\State \textbf{Initialization}: $n,~m,~\epsilon,~K,~k=0,~\bm{\lambda}^{0}=\bm{0}\in \mathbb{R}^m,~Index=0,~Iter=0$.
\For{$Iter\leq K$}
\State{$Iter=Iter+1$}.
\State{$k=k+1$}.
\State{$Index=0$.}
\If{ $k> m$}
\State{$k=1$}.
\EndIf
\State{Set $list(k)$ and generate the corresponding problem \eqref{SSQPQC}.  }
\State{Solve problem \eqref{SSQPQC} by Algorithm \ref{alg1} to obtain the multiplier $\lambda_{k}$ and $\bm{y}^{*}({\lambda}_{k})$. }
\State{Set $\lambda^{Iter}_{k}={\lambda}_{k}$ and  $\bm{y}^{*}=\bm{y}^{*}({\lambda}_{k})$.}
\For{$i =1,\dots,m$}
\If{$g_{i}(\bm{y}^{*}) \leq \epsilon~\text{or}~|\lambda_{i}g_{i}(\bm{y}^{*})|\leq\epsilon$}
\State{\textit{Index}= \textit{Index}+1}.
\EndIf
\EndFor
\If{$Index=m$}
\State{Break}.
\EndIf
\EndFor
\State \textbf{Output}: The optimal dual multiplier $\bm{\lambda}^{*}=\bm{\lambda}^{Iter}$ and the optimal solution $\bm{y}^{*}$ of problem \eqref{SQPQC}.
\end{algorithmic}
\end{algorithm}

\begin{proposition}
Applying  the Main Algorithm to solve problem \eqref{SQPQC} with $m$ quadratic constraints, if $Index = m$, then $\lambda_{k}$ and $\bm{y}^*(\lambda_{k})$ form an $\epsilon$-solution of the \eqref{KKT} system. 
\end{proposition}
\begin{proof}
After obtaining $\lambda_{k}$ and $\bm{y}^*(\lambda_{k})$,  $Index=m$ indicates that $|\lambda_{i}g_{i}(\bm{y}^*(\lambda_{k}))|\leq \epsilon$ and $g_{i}(\bm{y}^*(\lambda_{k}))\leq \epsilon $ for $i=1,...,m$. Combining with the fact that the primal-dual pair $(\lambda_{k},\bm{y}^*(\lambda_{k}))$ solves the corresponding system \eqref{SKKT}, it also solves the system \eqref{KKT}. 
\end{proof}




In the following, we will show that Main Algorithm is actually a dual coordinate ascent algorithm for problem \eqref{SQPQC} with the cyclic updating rule. To this end, we write the dual problem of \eqref{SQPQC} as follows. 
\begin{equation}\tag{DSQPQC}\label{DSQPQC}
    \max_{\bm{\lambda}\geq 0} ~ L(\bm{\lambda}),
\end{equation}
where $L(\bm{\lambda})=\min_{\bm{y}\in[\bm{l}, \bm{u}]}  \left\{f(\bm{y})+\sum_{i=1}^{m} \lambda_{i} g_{i}(\bm{y})\right\}$ is the Lagrangian function and $\lambda_i$ is the dual multiplier of constraint $g_{i}(\bm{y})\leq 0$ for $i=1,...,m$. The dual coordinate ascent algorithm for problem \eqref{SQPQC} optimizes the dual problem \eqref{DSQPQC} in one coordinate, say $\lambda_k$ for some $k\in \{1,....,m\}$, at a time while keeping other coordinates $\lambda_{i}$, $i\in list(k)$, fixed. The cyclic updating rule means that the algorithm cycles through the dual variables in sequence. 
\begin{proposition}
The solution $\lambda_{k}$ to the subsystem \eqref{SKKT} is the optimal solution of the subproblem in the dual coordinate ascent algorithm for problem \eqref{SQPQC} with the cyclic updating rule.
\end{proposition}
\begin{proof}
In each iteration, the dual coordinate ascent algorithm for problem \eqref{SQPQC} solves the following subproblem for a given $k\in \{1,...,m\}$:
\begin{align*}
   \max_{\lambda_k \ge 0} ~L(\lambda_{k}) := \min_{\bm{y}\in[\bm{l}, \bm{u}]} & \left\{ f(\bm{y})+\sum_{i\in list(k)} \hat{\lambda}_{i} g_{i}(\bm{y}) +\lambda_{k}g_{k}(\bm{y})\right\},
\end{align*}
 where $\hat{\lambda}_i$ for $i\in list(k)$ are fixed constants. Note that the KKT condition of the above subproblem is exactly the subsystem \eqref{SKKT}. Therefore, the solution $\lambda_{k}$ to  subsystem \eqref{SKKT}  solves the subproblem in the dual coordinate ascent algorithm for problem \eqref{SQPQC}. 
\end{proof}

Before showing the convergence proof of the Main Algorithm, we take care of the differentiability of $L(\bm{\lambda})$ in the next result. 
\begin{proposition}\label{prop4}
     $L(\bm{\lambda})$ is differentiable at any $\bm{\lambda}\geq 0$.
\end{proposition}
\begin{proof}
  { Under Assumption 1, for any given $\bm{\lambda}\ge 0$, the optimal solution $\bm{y}^{*}(\bm{\lambda})$ is unique. Since $f(\cdot)$ and $g_{i}(\cdot), ~i=1,\dots,m$, are continuous convex functions, it follows from Theorem 35.8 in \cite{rockafellar1997convex} that $D(\bm{\lambda})$ is continuously differentiable.}
\end{proof}
Now we are ready to give the main result of this section. 
\begin{theorem}
Under Assumption 1, Assumption 2 and the condition $$\left\{\bm{{y}} \mid g_{i}(\bm{{y}})=0\right\} \cap \{\bm{{y}} \mid {y}_{j}\in \{l_{j},u_{j},-\frac{\beta_{i}^{j}}{2\theta_{i}^{j}}\},~ j=1,\dots, n\}=\emptyset, ~ i=1,\dots,m,$$ for problem \eqref{SQPQC}, if $\{\bm{\lambda}^{Iter}\}$ is a sequence of solutions generated by Main Algorithm, then every limit point of $\{\bm{\lambda}^{Iter}\}$ is an optimal dual solution $\bm{\lambda^*}$ of problem \eqref{DSQPQC}, and $\bm{y}^*$ is an optimal solution of problem \eqref{SQPQC}.
\end{theorem}
\begin{proof}
Under the given conditions, it follows from Proposition \ref{Prop2.4} that the dual problem $\max_{\lambda_{k}\geq0 } L(\lambda_{k})$ in each iteration of Main Algorithm has a unique optimal solution. Furthermore, according to Proposition \ref{prop4}, the objective function $L(\bm{\lambda})$ is concave and differentiable over the convex and closed feasible domain of $\bm{\lambda}\ge0$.  Leveraging Proposition 6.5.1 in \cite{bertsekas2015convex}, we know that every limit point of $\{\bm{\lambda}^{Iter}\}$ is an optimal solution to the dual problem \eqref{DSQPQC}. Hence, the corresponding primal solution $\bm{y}^*$ is an optimal solution to problem \eqref{SQPQC}. 
\end{proof}

\section{Numerical experiments} \label{sec4}
In this section, we conduct computational experiments to test the efficiency of the proposed algorithm. The  algorithm is implemented in MATLAB R2019a on a laptop equipped with an Intel Core i5-8300H CPU with 8 GB RAM and Windows 10 OS. We compare Main Algorithm with the  commercial solver Gurobi (version 10). The settings of Gurobi are set to be the default values, except that the maximal runtime is limited to 7,200 cpu-seconds (2 hours). The tolerance error $\epsilon$ in  Main Algorithm  and Algorithm 1 is set to be $10^{-6}$,  which is same as the default feasibility  and optimality tolerance in Gurobi. The maximum number of iterations $K$ is set to be $1000$. To generate random instances of problem \eqref{SQPQC}, we modify the generator in \cite{Luo2023}. The implementation details of the generator are described as follows:
\begin{itemize}
	\item $\delta_j \in U[0,1]$ (i.e., uniformly distributed within the interval [0,1]) and $\bm{\alpha}_j \in U[-5,-2]$ for $j=1, \ldots, n$.
	
	\item $\theta^{i}_j \in U[0,2]$ for $i= 1,\dots,m$ and $j=1, \ldots, n$.  $\bm{\beta}^{i}_j \in U[0,5]$ for $i=1,\dots,m$  and  $j=1, \ldots, n$.
	
	\item $l_j=-1$ and $u_{j}=1$ for $j=1, \ldots, n$.
	
	\item To guarantee feasibility, we set $\sigma_{i} \in U[-v_{i}-1,-v_{i}]$ with $v_{i}=\bm{y_0}^T \Theta_{i} \bm{y_0}+\bm{\beta}_{i}^T \bm{y_0}$, where $\bm{y_0}$ is randomly drawn from $[-1,1]^{n}$.
\end{itemize}

To test the performance of the proposed algorithm for large-scale instances, we first fix the number of constraints $m$ to two and increase the size of the problem $n$ from one thousand to one million. For each problem size, we randomly generate ten instances. The numerical results are summarized in Table \ref{tab1}, in which the column “Solved” stands for the number of solved instances within the time limit of two hours, column “Time”  for the averaged computational time (in seconds) over the solved instances, column “Iter” for the average number of iterations, and column “Gap” for the average gap between the objective values returned by  Main Algorithm and Gurobi. 

\begin{table}[H]
	\caption{Comparison of   Main Algorithm with Gurobi for random instances with $m=2$.} \label{tab1}
	\centering
	\begin{tabular}{lccc|ccc}
		\toprule
		\multirow{2}{*}{$n$ }& \multicolumn{3}{@{}c@{}}{ Main Algorithm} & \multicolumn{2}{c}{Gurobi} & \multirow{2}{*}{Gap ($\times10^{-6}$)}\\
		\cmidrule{2-6}
		  & Solved    & Time   & Iter    & Solved  & Time             \\
		\midrule
		1000 &10 & 0.036 & 34  & \multicolumn{1}{|c}{10}& 0.133 &0.106\\ 
		5000 &10& 0.284 & 35  & \multicolumn{1}{|c}{10}& 2.414 &0.102\\ 
		10000 &10& 0.627 & 37  & \multicolumn{1}{|c}{10}& 8.168 &0.122 \\ 
		50000 &10& 2.610 & 39  & \multicolumn{1}{|c}{10}& 80.846 &0.102\\  
		100000 &10& 4.309 & 41  & \multicolumn{1}{|c}{10}& 316.070 &0.119\\ 
		500000 &10& 44.201 & 43  & \multicolumn{1}{|c}{0}&  ~-$^{1}$ & -$^{1}$\\ 
		1000000 &10& 93.048 & 44  & \multicolumn{1}{|c}{0}& ~-$^{1}$ &-$^{1}$\\ 
		\bottomrule
	\end{tabular}
	\begin{tablenotes}
		\footnotesize
		\item ${}^1$All instances failed solving by Gurobi within the time limit of two hours.
	\end{tablenotes}
\end{table}

The results in  Table \ref{tab1} clearly indicate that  the proposed  algorithm outperforms Gurobi solver by at least one order of magnitude in terms of running time. 
The average difference in objective values between the proposed algorithm and Gurobi is less than $10^{-6}$, indicating that  Main Algorithm indeed returns an optimal solution of problem \eqref{SQPQC}.

As the size of the problem increases, the superiority of  the proposed Main Algorithm becomes more apparant. When the problem size is over half a million, Gurobi fails to solve any instances within the two-hour time limit, whereas the proposed algorithm can solve all instances within 100 seconds. This result shows that the proposed algorithm is  promising  for solving large-scale \eqref{SQPQC} problems. 

Next, we assess the sensitivity of the proposed algorithm to the number of constraints. We vary the number of quadratic constraints $m$ from $2$ to $20$ while fixing problem size $n=50000$. The numerical results of the proposed algorithm and Gurobi are summarized in Table \ref{tab2}. The results show that the average solution time of Main Algorithm and Gurobi increases as $m$ becomes larger. However, the proposed algorithm still runs significantly faster than Gurobi by at least one order of magnitude. Gurobi does not solve some instances when $m = 10$, while the proposed algorithm begins to fail from $m=20$.  The numerical results show that the proposed algorithm is efficient for problem \eqref{SQPQC} with a few quadratic constraints and boxed variables, but may not be suitable for problems with many constraints. 


\begin{table}[htbp!]
	\caption{Comparison of  Main Algorithm with Gurobi for random instances with $n = 50000$ and different   numbers of constraints $m$.} \label{tab2}
	\centering
	\begin{tabular}{cccc|ccc}
		\toprule
	\multirow{2}{*}{$m$ }	& \multicolumn{3}{@{}c@{}}{ Main Algorithm} & \multicolumn{2}{c}{Gurobi} & \multirow{2}{*}{Gap ($\times10^{-6}$)}\\
		\cmidrule{2-6}
		     & Solved & Time    & Iter     & Solved  & Time           \\
		\midrule
	2 &10& 2.610 & 39 & \multicolumn{1}{|c}{10} & 80.846  & 0.102\\ 
	3 &10& 3.551 & 40 & \multicolumn{1}{|c}{10} & 121.900  &0.165\\ 
	4 &10& 6.791 & 55 & \multicolumn{1}{|c}{10}  & 157.284 & 0.093\\ 
	5 &10& 11.331 & 72 & \multicolumn{1}{|c}{10} & 200.035  & 0.108\\ 
	6 &10& 18.716 & 91 & \multicolumn{1}{|c}{10} & 291.137  & 0.048\\ 
	7 &10& 23.815 & 101 & \multicolumn{1}{|c}{10} & 387.202  & 0.126\\ 
    10 & 10 & 86.028 & 215  &\multicolumn{1}{|c}{3} & 496.504&0.146 \\ 
    15 & 10 & 274.298 & 325 &\multicolumn{1}{|c}{0} &- &-\\ 
    20 & 3 & 485.582 & 496 &\multicolumn{1}{|c}{0} &- &-\\     
		\bottomrule
	\end{tabular}
\end{table}

\section{Conclusion}\label{sec5}
In this paper, we   propose an efficient solution method for solving the  convex separable quadratic optimization problem \eqref{SQPQC} in large size. Utilizing an efficient algorithm for solving convex separable quadratic problem with one quadratic constraint and boxed variables, we have developed a dual coordinate ascent algorithm for solving problem \eqref{SQPQC} via an iterative resolution scheme for the KKT system, and provided a convergence proof of the proposed algorithm. The superior performance of the proposed algorithm for solving large-size instances of problem \eqref{SQPQC} with a few quadratic constraints and boxed variables has been illustrated by  comparisons with the widely used commercial solver Gurobi. 

There are two possible directions for  future study. One is to extend the proposed algorithm for solving  nonconvex separable quadratic programs with multiple quadratic constraints. The other one is to develop efficient algorithms for solving large-scale convex quadratic inseparable problems with special structures.

\section*{Data Availability}

The data that support the findings of this study are available from the corresponding author upon request.

\section*{Acknowledgement}
Lu's research has been supported by the National Natural Science Foundation of China Grant No. 12171151. Deng's research has been supported by the National Natural Science Foundation of China Grant No. T2293774, by the Fundamental Research Funds for the Central Universities E2ET0808X2, and by the grant from MOE Social Science Laboratory of Digital Economic Forecast and Policy Simulation at UCAS.

    



\begin{thebibliography}{99}





\bibitem{bertsekas2015convex} D.P. Bertsekas, \emph{Convex Optimization Algorithms}, Athena Scientific, Nashua, 2015.

\bibitem{bitran1981disaggregation} G.R. Bitran and A.C. Hax, Disaggregation and resource allocation using convex knapsack problems with bounded variables, \emph{Manage. Sci.}  {27} (1981)  431--441.


\bibitem{Chen2014} J. Chen, L. Feng, J. Peng and Y. Ye, Analytical results and efficient algorithm for optimal portfolio deleveraging with market impact, \emph{Oper. Res.}  {62} (2014)  195--206.

\bibitem{cominetti2014newton} R. Cominetti, W.F. Mascarenhas and P.J.S. Silva, A Newton’s method for the continuous quadratic knapsack problem, \emph{Math. Program. Comput.}  {6} (2014) 151--169.

\bibitem{dai2006new} Y. Dai and R. Fletcher, New algorithms for singly linearly constrained quadratic programs subject to lower and upper bounds, \emph{Math. Program.}  {106} (2006) 403--421.

\bibitem{pjo1} Y. Dai and K. Schittkowski, A sequential quadratic programming algorithm with non-monotone line search, \emph{Pac. J. Optim.} 4 (2008) 335-358.



\bibitem{edirisinghe2016efficient} C. Edirisinghe and J. Jeong, An efficient global algorithm for a class of indefinite separable quadratic programs, \emph{Math. Program.}  {158} (2016) 143--173.

\bibitem{edirisinghe2017tight} C. Edirisinghe and J. Jeong, Tight bounds on indefinite separable singly-constrained quadratic programs in linear-time, \emph{Math. Program.}  {164} (2017) 193--227.

\bibitem{fukushima2003sequential} M. Fukushima, Z.Q. Luo and P. Tseng, A sequential quadratically constrained quadratic programming method for differentiable convex minimization, \emph{SIAM J. Optim.}  {13} (2003) 1098--1119.

\bibitem{hochbaum1995strongly} D.S. Hochbaum and S.P. Hong, About strongly polynomial time algorithms for quadratic optimization over submodular constraints, \emph{Math. Program.}  {69} (1995)  269--309.

\bibitem{pjo2} A. Kato, N. Yasushi and Y. Hiroshi, Global and superlinear convergence of inexact sequential quadratically constrained quadratic programming method for convex programming, \emph{	Pac. J. Optim.} 8 (2012) 609-629.

\bibitem{li2023efficient} S. Li, Z. Deng, C. Lu, J. Wu, J. Dai and Q. Wang, An efficient global algorithm for indefinite separable quadratic knapsack problems with box constraints, \emph{Comput. Optim. Appl.}  {86} (2023) 241--273.


\bibitem{lobo1998applications} M.S. Lobo, L. Vandenberghe, S. Boyd and H. Lebret, Applications of second-order cone programming, \emph{Lin. Alg. Appl.}  {284} (1998) 193--228.

\bibitem{Luo2023} H. Luo, X. Zhang, H. Wu and W. Xu, Effective algorithms for separable nonconvex quadratic programming with one quadratic and box constraints, \emph{Comput. Optim. Appl.} {86} (2023) 199--240. 

\bibitem{megiddo1984linear} N. Megiddo, Linear programming in linear time when the dimension is fixed, \emph{J. ACM}  {31} (1984) 114--127.

\bibitem{megiddo1993linear} N. Megiddo and A. Tamir, Linear time algorithms for some separable quadratic programming problems, \emph{Oper. Res. Lett.}  {13} (1993),  203--211.

\bibitem{mehrotra1991method} S. Mehrotra and J. Sun, A method of analytic centers for quadratically constrained convex quadratic programs, \emph{SIAM J. Numer. Anal.}  {28} (1991)  529--544.



\bibitem{monteiro2000polynomial} R.D.C. Monteiro and T. Tsuchiya, Polynomial convergence of primal-dual algorithms for the second-order cone program based on the MZ-family of directions, \emph{Math. Program.}  {88} (2000) 61--83.


\bibitem{rockafellar1997convex} R.T. Rockafellar, \emph{Convex Analysis}, Princeton University Press, Princeton, 1997.

\bibitem{shetty1990parallel} B. Shetty and R. Muthukrishnan, A parallel projection for the multicommodity network model, \emph{J. Oper. Res. Soc.}  {41} (1990) 837--842.

\bibitem{tseng1990dual} P. Tseng, Dual ascent methods for problems with strictly convex costs and linear constraints: A unified approach, \emph{SIAM J. Control Optim.}  {28} (1990) 214--242.

\bibitem{tsuchiya1999convergence} T. Tsuchiya, A convergence analysis of the scaling-invariant primal-dual path-following algorithms for second-order cone programming, \emph{Optim. Methods Softw.}  {11} (1999) 141--182.



































\end{thebibliography}

\hbox to14cm{\hrulefill}\par
\vspace{2mm}
\noindent \textsc{Shaoze Li}\\
School of Economics and Management\\
University of Chinese Academy of Sciences, Beijing 100190, China\\
E-mail address: shaoze\_li@163.com \\
\\
\textsc{Junhao Wu}\\
School of Economics and Management\\
North China Electric Power University, Beijing 102206, China \\
E-mail address: junhao\_wu@163.com \\
\\
\textsc{Cheng Lu}\\
School of Economics and Management\\
North China Electric Power University, Beijing 102206, China \\
E-mail address: lucheng1983@163.com \\
\\
\textsc{Zhibin Deng (corresponding author)}\\
School of Economics and Management \\
University of Chinese Academy of Sciences, Beijing 100190, China\\
MOE  Social Science Laboratory of Digital Economic Forecasts and Policy Simulation at UCAS, Beijing, 100190, China\\
E-mail address: zhibindeng@ucas.edu.cn \\
\\
\textsc{Shu-Cherng Fang}\\
Department of Industrial and Systems Engineering \\
North Carolina State University, Raleigh 27695-7906, USA \\
E-mail address: fang@ncsu.edu \\



\end{document}